\documentclass[10pt]{article}
\usepackage{mathrsfs}
\renewcommand{\baselinestretch}{1}

\newtheorem{Th}{Theorem}[section]
\newtheorem{Le}{Lemma}[section]
\newtheorem{Co}{Corollary}[section]
\newtheorem{De}{Definition}[section]
\newtheorem{Pro}{Proposition}[section]
\newtheorem{Rem}{Remark}[section]
\newtheorem{Ex}{Example}[section]
\newcommand{\bth}{\begin{Th}}
\newcommand{\eth}{\end{Th}}
\newcommand{\ble}{\begin{Le}}
\newcommand{\ele}{\end{Le}}
\newcommand{\bco}{\begin{Co}}
\newcommand{\eco}{\end{Co}}
\newcommand{\bde}{\begin{De}}
\newcommand{\ede}{\end{De}}
\newcommand{\bpr}{\begin{Pro}}
\newcommand{\epr}{\end{Pro}}
\newcommand{\bre}{\begin{Rem}}
\newcommand{\ere}{\end{Rem}}
\newcommand{\bex}{\begin{Ex}}
\newcommand{\eex}{\end{Ex}}
\newcommand{\beq}{\begin{equation}}
\newcommand{\eeq}{\end{equation}}
\newcommand{\beqn}{\begin{eqnarray}}
\newcommand{\eeqn}{\end{eqnarray}}
\newcommand{\be}{\begin{eqnarray*}}
\newcommand{\ee}{\end{eqnarray*}}

\def\mysection#1{\setcounter{equation}{0}
{\renewcommand{\thesection}{\arabic{section}}\section{#1}}
\renewcommand{\thesection}{\arabic{section}}}
\usepackage{cite}
\begin{document}
\title{
{\renewcommand{\baselinestretch}{2}\Huge \bf ON BOUNDARY PROBLEMS FOR REGULAR FUNCTIONS IN
HYPERCOMPLEX ANALYSIS}
 \thanks{This work is
supported by  Macao Science and Technology Development Fund, MSAR. Ref. 041/2012/A,
 NNSF of China (\#10471107) and SRFDP of Higher Eduction of China (\#20060486001).}}
\author{ Jinyuan Du$^{1}$ \, \& \,  Pei Dang$^{\,2}$\\
\small 1. Department of Mathematics, Wuhan University, Wuhan
$430072$, P. R. China\\
\small  E--mail: \hspace{2mm}jydu@whu.edu.cn\\
\small 2.  Faculty of Innovation Engineering, Macau University of Science and Technology, Macao, China\\
\small  E--mail: \hspace{2mm}pdang@must.edu.mo}
\date{}
\maketitle



\begin{minipage}{150mm}
\abstract{ \,\,In this article, the authors survey and review the studies of boundary value
problems for regular functions in Clifford analysis, which include
theoretical foundations and useful methods. Its theoretical bases
consist of the generalized Cauchy theorem, the  generalized Cauchy
integral formula, the Painlev\'{e} theorem and boundary behaviors of
the Cauchy type integrals, as well as various integral
representations. Certain boundary value problems in the Clifford
algebra setting and singular integral equations are introduced.
}\vspace{3mm}

\noindent {\bf Keywords.} {\small  Cauchy-Goursat theorem,  Painlev\'{e} theorem,  Boundary behavior,
Riemann boundary value problem, Dirichlet boundary value problem, Singular integral equation.\vspace{1mm}}

\noindent{{\bf 2000 MR Subject Classification:} \,\small 30G35;
31B05; 35C10.}
\end{minipage}


\mysection{Introduction}

It has been proved that the function theory over Clifford algebra is
an appropriate setting to generalize many aspects of function theory
of one complex variable to higher dimensions
\cite{BDS-RNM-82,DelangheSS-92,GurlebeckS-90,GurlebeckS-97,w91}.
Boundary value problems (BVPs) for analytic functions in the
classical complex analysis is an important branch of mathematical
analysis \cite{Lu1,Mu1,Ga1}, and has been actively studied due to
its theoretical elegance and widespread applications in physics and
other subjects such as elastic theory, hydromechanics,  fracture
mechanics, etc. \cite{lu2,m2,Vekua}.

 It is natural that we hope to develop corresponding theories
on boundary value problems for regular functions in the hypercomplex
analysis similarly to those for analytic functions in the classical
complex analysis. In fact, some results on boundary value problems
for analytic functions in the classical complex analysis have been
generalized to those for regular functions in Clifford analysis
\cite{Xu87,Xc87,Xu90,Xu92,Xz93,Bl00,Bl01,S91,S95,SV95, Stern,Be96,
Be962,Be98,h1,qy1,Gz10,Duzh-01,DuGong-04,gd,jd12,sd2,sd3,bd }.
However, to our knowledge, the foundation for boundary problems in
Clifford analysis has not been established completely. For example,
the generalized Cauchy theorem (i.e., Cauchy-Goursat theorem) and
 generalized Cauchy integral formula had not been conscientiously
treated in Clifford analysis for a long time after the Cauchy
theorem and Cauchy integral formula were obtained in Clifford
analysis \cite{BDS-RNM-82, dz2+, bdz }.
 The generalized Cauchy theorem and the
 generalized Cauchy integral formula play crucial
roles in the theory of boundary value problems for analytic
functions in the classical complex analysis, so do their analogues
in Clifford analysis, still called Cauchy-Goursat theorem and the
 generalized Cauchy integral formula.
Due to the necessity of Cauchy-Goursat theorem  and  the generalized
Cauchy integral formula, many studies tacitly assumed their
existence. Before they were strictly proved, however, Painlev\'{e} theorem is
lack of rigorous foundation.  In our opinion, the works on BVPs for
regular functions lose theoretical support before these basic
theories (e.g., Painlev\'{e} theorem) are strictly set up.

Our aim here is to bring together some dispersed materials, not to
give a complete survey. We will point out some major contributions
on the basic theory related to boundary value problems for regular
functions and the recent results and methods to solve boundary value
problems for regular functions. In next section, we introduce some
major contributions on boundary behavior of the Cauchy type
integral. In \S 3, we point out  the recent development on boundary
value problems for regular functions. Then, the jump Riemann
boundary value  problems for regular functions in Clifford analysis
are formulated and some open problems for general Riemann boundary
value problems are given in \S 4. In  \S 5, we discuss some singular
integral equations and the inversion formula for Cauchy principal
value. In \S 6,  we establish the generalized Sochocki--Plemelj
formula  and solve the Dirichlet problems by using it.

Throughout this paper, all notations and terminologies are referred
to \cite{du+}. \vspace{2mm}

\mysection{Boundary behavior of Cauchy type integrals}

Let $C(V_{n})$ be the Clifford algebra obtained by algebraically
spanning the linear subspace ($1$-vector space) $V_{n}\!=\!\mbox{\rm
span}\{e_1,\cdots,e_n\}$ \cite{BDS-RNM-82,du+}, of which  basis
are written as $\left\{ e_A, A\!=\!(h_1,\cdots, h_r)\!\in\!\mbox{$\cal
P$}\!N, 1\leq h_1<\cdots<h_r \leq n \right\},$ where $N$ stands for
the set $\{1,\cdots, n \}$ and $\mbox{$\cal P$}\!N$ denotes the
family of all order-preserving subsets of $N$ in the above fixed
way. We assume that $\Omega^+$ is a bounded open domain in
$\mathcal{R}^{n+1}$ with a $C^1$-smooth boundary
$\Gamma=\partial\Omega^+,$ that
 is oriented by the exterior normal, and
 $\Omega^-=\mathcal{R}^{n+1}\backslash \{\Omega^+\cup\Gamma\}$.
 We  agree  that
  the regular (entire)  function  always names one of the left regular (entire) function and the right regular (entire) function,
 unless otherwise stated, we often abbreviate $\lq\lq$ left regular (entire) function"
or $\lq\lq$right regular (entire) function" as $\lq\lq$ regular
(entire) function". \vspace{1.5mm}

\bde \label{sectionally}
A function $F$ is said to be sectionally regular with $\Gamma$ as its jump surface, if it is  regular in $\Omega^{+}$ and
$\Omega^{-}$, and has positive and negative boundary values when $w$ tends to any point $x\in \Gamma$ from $\Omega^{+}$ and
$\Omega^{-}$, i.e., there exist

\beq  \label{21}
F^{+}(x)=\lim_{\Omega^{+} \ni w\rightarrow x}F(w),\,\,\,\,\,\,\,\,\,F^{-}(x)=\lim_{\Omega^{-} \ni w\rightarrow x}F(w).
\eeq
\ede
\bre\label{R21}
{\rm It is easy to prove that $F^{+}\in C\!\left(\overline{\Omega^{+}},C(V_n)\right)$ and
$F^{-}\in C\!\left(\overline{ \Omega^{-}},C(V_n)\right)$.
}\ere

Now we introduce some results related to the jump Riemann boundary value problem.\vspace{1mm}

{\bf  Jump Riemann boundary value problem} (jump problem).
Find a sectionally  regular function $\Phi$,  with $\Gamma$ as its jump surface,  satisfying the boundary value condition
\beq \label{22}
  \Phi^+(x)=\Phi^-(x)+g(x),\,\,\,\, x\in\Gamma,\eeq
where $g\in H(\Gamma)$ (i.e., H\"{o}lder continuous on $\Gamma$, seen below in Detail).\vspace{1mm}

To find a special solution of the jump problem (\ref{22}), we
introduce the Cauchy type integrals with the density $f$ as follows.
\beq \label{23} \Big(C[f]\Big)(w)=
 \displaystyle\frac{1}{\bigvee_{n+1}}\int_{\Gamma} E(x-w)\,\mathrm{d}\sigma f(x),\,\,\,w\in\mathcal{R}^{n+1}\backslash\Gamma
   \eeq
   and
      \beq\label{24}
\Big([f]C\Big)(w)=\displaystyle\frac{1}{\bigvee_{n+1}}\int_{\Gamma}f(x)\, \mathrm{d}\sigma\, E(x-w), \,\,\,w\in\mathcal{R}^{n+1}\backslash\Gamma,
   \eeq
   where $E$ is the Cauchy kernel function, given by
  \beq\label{25}
    E(x-w)=\frac{\overline{x-w}}{|x-w|^{n+1}},\,\,\,\,\,\,x\not=w,
  \eeq

   \beq\label{26}
{\bigvee}_{n+1}=\displaystyle\frac{2 \pi^{(n+1)/2}}{\Gamma((n+1)/2)}\vspace{-1mm}
\eeq
  is the area of the unit sphere
  \beq\label{27}
  S^{n}=\Big\{x,\,\, |x|=1,\, x=\left(x_0,x_1,\cdots,x_n\right)\in\mathcal{R}^{n+1} \Big\},
\eeq
and
\beq\label{28}
 \mbox{\rm d}\sigma=\sum_{k=0}^{n}(-1)^{k}e_k\mbox{\rm d}\widehat{x}_{k}=
 \sum_{k=0}^{n}(-1)^{k}e_k\mbox{\rm d}x_{0}\wedge\cdots\wedge\mbox{\rm
d}x_{k-1}
 \wedge\mbox{\rm d}x_{k+1}\wedge\cdots\wedge\mbox{\rm d}x_n.
 \eeq
\bth [Theorem for the Cauchy type integral {\rm\cite{du+,dx,dxz}}]
\label{c1} If $f$ is absolutely integrable on $\Gamma$, then the
Cauchy type integrals $C[f]$ and $[f]C$
 are, respectively, left and right regular. Moreover,
\beq \label{29}
 \limsup\limits_{|x|\rightarrow +\infty}|x|^{m}\left|\big(C[f]\big)(x)\right|<+\infty, \,\,\,\,\
 \,\,\,\,\limsup\limits_{|x|\rightarrow +\infty}|x|^{m}\left|\big([f]C\big)(x)\right|<+\infty,
 \eeq
 in particular,
 \beq \label{210}
 \lim_{|x|\rightarrow +\infty}|x|^{m-1}\big|(C[f])(x)\big|=0,\,\,\,\,
 \,\,\,\,\lim_{|x|\rightarrow +\infty}|x|^{m-1}\big|([f]C)(x)\big|=0.
 \eeq
\eth

 Theorem \ref{c1} is based on the following rule.

 \bth[Leibniz rule {\rm\cite{dx,dxz}}] Suppose that $\Omega$ is an open set in $\mbox{$\cal R$}^{n+1}$ and $\phi
 \in C(\Gamma\times\Omega)$. Let
\beq \label{211}
\Phi(w)=\int_{\Gamma}\phi(x,w)\mbox{\rm d}\sigma f(x),\,\,\,\,
w\in\Omega,
\eeq
if  $D[\phi]\in C(\Gamma\times\Omega)),$
then
\beq \label{212}
\big(D[\Phi]\big)(w)=\displaystyle\int_{\Gamma}\big(D[\phi]\big)(x,w)\mbox{\rm
d}\sigma f(x),\,\,\,\,
w\in\Omega.
\eeq
Let
\beq \label{213}
\Psi(w)=\int_{\Gamma}f(x) \mbox{\rm d}\sigma \psi(x,w),\,\,\,\,
w\in\Omega,
\eeq
if  $[\phi]D\in C(\Gamma\times\Omega)),$
then
\beq \label{214}
\big([\Psi]D\big)(w)=\displaystyle\int_{\Gamma}f(x) \mbox{\rm
d}\sigma \big([\psi]D\big)(x,w),\,\,\,\,
w\in\Omega.
\eeq
\eth

 Now we dress slightly up the Cauchy type integral.
 After that, we will see that it is a copy of the Cauchy type
 integral in one complex variable.

 First we introduce a mapping
 \beq\label{215}
 {\rm capital}_{n}\!:\,x=\big(x_0,x_1,\cdots,x_n\big)
\longmapsto \displaystyle\sum\limits_{i=0}^{n}x_ie_i=X.
\eeq
Obviously, it is
a proper isomorphism between $\mbox{$\cal R$}^{n+1}$ and the linear
subspace $\mbox{\rm span}\{e_0,e_1,e_2,\cdots,e_n\}$ of
$C\left(V_{n}\right),$ in which the elements $\lambda=\lambda_0+\lambda_1e_1+\cdots+\lambda_n e_n$ with $\lambda_0,\cdots,\lambda_n\in\mathcal{R}$ are called paravectors.
  In this view,
$E(x)=
\left|x\right|^{1-n}
X^{-1},$ where
$X^{-1}$ is the inverse element of
$ X.$ Since there is no multiplication in
$\mathcal{R}^{n+1},$ it is not confused that
$x^{-1}$ is treated as
$X^{-1}$. Moreover, if
$\alpha,\beta\in C\left(V_{n}\right)$ and
 there exists the inverse $\alpha^{-1}$,
we denote
\beq\label{216}
\alpha^{-1}\beta=
\mbox{$|$}\hspace{-1mm}\displaystyle\frac{\mbox{}\beta
\mbox{}}{\mbox{}\,\,\,\alpha\,\,\,\mbox{}},\,\,\,\,\,\,\,\, \beta
\alpha^{-\!1}=\displaystyle\frac{\mbox{}
\beta\mbox{}}{\mbox{}\,\,\,\alpha\,\,\,\mbox{}}\hspace{-1mm}\mbox{$|$},
\eeq
 which are, respectively, called
left division and right division.
In this way (\ref{23}) and (\ref{24}) may be rewritten as
\beq
\big(C[f]\big)(w)=\frac{1}{\bigvee_{n+1}}\int_{\Gamma}
\mbox{$|$}\hspace{-1mm}\displaystyle\frac{\mbox{}\,\mbox{\rm d}\sigma\,f(x)\mbox{}}
{\mbox{}\,\left|x-w\right|^{n-1}
(x-w)\,\mbox{}}
, \hspace{4mm}\,\,\,w\in\mathcal{R}^{n+1}\backslash\Gamma, \eeq
and
\beq
([f]C)(w)=\displaystyle\frac{1}{\bigvee_{n+1}}\int_{\Gamma}
\displaystyle\frac{\mbox{}f(x)\,\mbox{\rm d}\sigma\,\mbox{}}
{\mbox{}\,\left|x-w\right|^{n-1}
(x-w)\mbox{}\hspace{1mm}}\hspace{-1mm}\mbox{$|$}, \hspace{4mm}\,\,\,w\in\mathcal{R}^{n+1}\backslash\Gamma.
 \eeq

Before discussing boundary behaviors of the Cauchy type integrals,
we need introduce the Cauchy principal value integral. Sketching a
sphere $O(t,\delta)$ with center at $t\in\Gamma$ and  sufficiently
small radius $\delta>0,$ we obtain a piece of surface
$\Gamma(t,\delta)$ from $\Gamma.$ Next we define the Cauchy
principal value integral by the limit \beq \label{219}
\Big(C[f]\Big)(t)=\displaystyle\frac{1}{\bigvee_{n+1}}
\displaystyle\int_{\Gamma} E(x-t) \mbox{\rm d}\sigma
f(x)=:\lim\limits_{\delta\rightarrow0^+}
\displaystyle\frac{1}{\bigvee_{n+1}}\,
\displaystyle\int_{\Gamma\setminus\Gamma(t,\,\delta)} E(x-t)
\mbox{\rm d}\sigma f\!\left(x\right),\hspace{4mm} t\in\Gamma. \eeq
Similarly, \beq \label{220} \Big([f]C\Big)(t)=
\displaystyle\frac{1}{\bigvee_{n+1}}\, \displaystyle\int_{\Gamma}
f(x) \mbox{\rm d}\sigma\,E(x-t)=:\lim\limits_{\delta\rightarrow0^+}
\displaystyle\frac{1}{\bigvee_{n+1}}\,
\displaystyle\int_{\Gamma\setminus\Gamma(t,\,\delta)} f(x) \mbox{\rm
d}\sigma\,E(x-t) ,\hspace{4mm} t\in\Gamma. \eeq \bex{\rm (see
\cite{dxz}) If $f\equiv1$, then the Cauchy principal value integral
\beq \label{221} \Big(C[1]\Big)(t)= \Big([1]C\Big)(t)=\frac{1}{2},
\hspace{6mm}t\in\Gamma.
 \eeq
 }
\eex

In general, the H\"{o}lder condition of the density function is a
sufficient condition for existence of the
 Cauchy principal value integrals.

\bde Let $f:D\rightarrow C(V_{n,s})$. If there exists a constant
$M>0$ such that, for any $x,y\in D$,
\beq \label{222}
\Big|f(x)-f(y)\Big|\leq M\Big|x-y\Big|^\mu \, \,\,\,\,\,(0<\mu\leq1),
\eeq
then we say f to be H\"{o}lder continuous on $D$, denoted by
$f\in H^{\mu}(D)$. The constant $M$ is called a H\"{o}lder
constant of $f$, $\mu$ is called H\"{o}lder index of $f$. If the order $\mu$ is not emphasized,
it may be written briefly by $f\in H(D)$ or $H$.
\ede

\bth [Existence of the Cauchy
principal integral {\rm \cite{dx,dxz}}] If $f\!\in\!H(\Gamma)$, then Cauchy
principal integrals $\big(C[f]\big)(t)$ and $\big([f]C\big)(t)$ exist for each
$t\in\Gamma$. More precisely,
\beq\label{223}
\Big(C[f]\Big)(t)=\Big(\Delta[f]\Big)(t)+\frac{1}{2}f(t),\,\,\,\,\,\,\,\,\Big([f]C\Big)(t)=\Big([f]\Delta\Big)(t)+\frac{1}{2}f(t),
\eeq
where
\beq \label{224}
\Big(\Delta[f]\Big)(t)=\displaystyle\frac{1}{\bigvee_{n+1}}\displaystyle\int_{\Gamma}\,E(x-t)\,
 \mbox{\rm
d}\theta\Big[f(x)
-f(t)\Big], \,\,\,\,\,
\Big([f]\Delta\Big)(t)=\displaystyle\frac{1}{\bigvee_{n+1}}\displaystyle\int_{\Gamma}\,
\,\Big[f(x)
-f(t)\Big] \mbox{\rm
d}\theta\,E(x-t)
\eeq
are two improper integrals.
 \eth

In earlier studies of boundary behavior of the Cauchy type
integrals, researchers
 always assumed that $\Gamma$ was a Liapunov surface \cite{HQW, GB,iv1}.
The hypothesis of the Liapunov surface is too strong for discussing the boundary behavior.
It is well known that
the condition of smooth curve is enough for the classical one complex variable case. The reason for
such difference comes from the common proof used to discuss the Clifford analysis case.
In fact, the discussion for the Clifford analysis case is similar to that for the classical one complex
 variable case. We have found the substitute of the chord-arc inequality so that the
 convergence of improper integrals in (\ref{224}) are easily proved.
 That is just the  concept of standard surface proposed by Du in \cite{dxz,dx}.

The Cauchy type integrals  (\ref{23}) and (\ref{24}) with the H\"{o}lder continuous density $f$
are just  sectionally  regular with $\Gamma$ as the jump surface. Their boundary values satisfy
the Sochocki--Plemelj formula and
Privalov--Muskhelishvili theorem.
We state these results in an unified framework,
which is collectively called $2$P theorem.
Introduce the following singular integral operators
\beq
\Big(S^{+}[f]\Big)(w)=\left\{
\begin{array}{lll}
\Big(C[f]\Big)(z),\hspace{1mm} &\mbox{\rm if}& w\in \Omega^{+},\\[4mm]
\displaystyle\frac{1}{2}f(t)+\Big(C[f]\Big)(t), \hspace{2.5mm}&\mbox{\rm if}& w=t\in\Gamma,
\end{array}
\right.
\eeq
\beq
\Big(S^{-}[f]\Big)(w)=\left\{
\begin{array}{lll}
\Big(C[f]\Big)(w), &\mbox{\rm if}& w\in \Omega^{-},\\[3mm]
-\displaystyle\frac{1}{2}f(t)+\Big(C[f]\Big)(t), &\mbox{\rm if}& w=t\in\Gamma,
\end{array}
\right.
\eeq

\beq
\Big([f]S^{+}\Big)(w)=\left\{
\begin{array}{lll}
\Big([f]C\Big)(w), &\mbox{\rm if}& w\in \Omega^{+},\\[3mm]
\displaystyle\frac{1}{2}f(t)+\Big([f]C\Big)(t), \hspace{2.5mm}&\mbox{\rm if}& w=t\in\Gamma,
\end{array}
\right.
\eeq

\beq
\Big([f]S^{-}\Big)(w)=\left\{
\begin{array}{lll}
\Big([f]C\Big)(w), &\mbox{\rm if}& w\in \Omega^{-},\\[3mm]
-\displaystyle\frac{1}{2}f(t)+\Big([f]C\Big)(t), &\mbox{\rm if}& w=t\in\Gamma.
\end{array}
\right.
\eeq

\bth[2P theorem {\rm \cite{dxz,dx}}]\label{2p}
If $f\in H^{\mu}(\Gamma)$, then
\beq
S^{+}[f], \,\,[f]S^+\in\left\{\begin{array}{lll}
H^{\mu}\Big(\overline{\Omega^+}\Big), &\mbox{\it if}& 0<\mu<1,\\[3mm]
H^{\nu}\Big(\overline{\Omega^+}\Big), &\mbox{\it if}& \mu=1,\,\, \nu<1,
\end{array}
\right.
\eeq
\beq
S^{-}[f], \,\,[f]S^-\in\left\{\begin{array}{lll}
H^{\mu}\Big(\overline{\Omega^-}\Big), &\mbox{\it if}& 0<\mu<1,\\[3mm]
H^{\nu}\Big(\overline{\Omega^-}\Big), &\mbox{\it if}& \mu=1,\,\, \nu<1.
\end{array}
\right.
\eeq
\eth

The 2P Theorem is a fundamental theorem of boundary behavior of
Cauchy type integrals.

\bre {\rm Noting the continuity of $S^{\pm}[f]$ and $[f]S^{\pm}$, at
once the Sochocki--Plemelj formulae for the boundary values of the
Cauchy type integrals $S[f]$ and $[f]S$ become a corollary of the 2P
Theorem.} \ere \bco [Sochocki--Plemelj formulae]\label{SP}
 If $f\in H(\Gamma)$, then, for $t\in\Gamma$,
\beq \label{231}
\Big(S[f]\Big)^{\!+}(t)
=\displaystyle\frac{1}{2}f(t)+\Big(C[f]\Big)(t),\hspace{5mm}
\Big(S[f]\Big)^{\!-}(t)
=-\displaystyle\frac{1}{2}f(t)+\Big(C[f]\Big)(t),
\eeq
\beq\label{232}
\Big([f]S\Big)^{\!+}(t)
=\displaystyle\frac{1}{2}f(t)+\Big([f]C\Big)(t),\hspace{5mm}
\Big([f]S\Big)^{\!-}(t)
=-\displaystyle\frac{1}{2}f(t)+\Big([f]C\Big)(t).
\eeq
\eco

\bre{\rm Theorem \ref{c1} with the Sochocki--Plemelj formulae (\ref{231}) and (\ref{232}) straightforward show that
$S[g]$ and $[g]S$ are, respectively, the special left regular solution and the special right regular solution of the jump problem (\ref{22}).
}
\ere

 The earliest result for the Sochocki--Plemelj formula
in Clifford analysis is given by Iftimie \cite{iv1}. His results took the Cauchy type integrals
over a Liapunov surface and the limits along non-tangent direction.
He followed the classical method in \cite{Mu1} to prove those results, that seems a little complicated. In fact, in 1980, Jinyuan Du has given a simple proof
of the Sochocki--Plemelj formulae, the Privalov theorem as well as Muskhelishvili
theorem for Cauchy type integrals in the classical one complex
variable case \cite{Lu1, du++}. The method in \cite{du++} still works  to prove
the Plemelj formulae and the Privalov--Muskhelishvili theorems for Cauchy type integrals in
 Clifford analysis case \cite{dxz}, even in the case of Cauchy type integrals with value in
 a universal Clifford algebral \cite{dx}. The 2P Theorem here
also improves the corresponding result in \cite{Z1}. In \cite{Z1},
Zhang obtained the result with the H\"{o}lder index $(\mu-\epsilon)$
using the method in \cite{du++}.

\bre{\rm
The Cauchy principle value integrals (\ref{219}) and (\ref{220}) are H\"{o}lder continuous  by 2P theorem.
}
\ere

\bco[H\"{o}lder continuity of the Cauchy principal integrals] If $f\!\in\!H(\Gamma)$, then $C[f]\!\in\!H(\Gamma)$ and $[f]C\in H(\Gamma)$.
In particular, if $f\in H^{\mu}(\Gamma),$ then $C[f]\in H^{\mu}(\Gamma)$ with $\mu<1;$ if $f\in H^{1}(\Gamma),$
then $C[f]\in H^{\nu}(\Gamma)$ for any $\nu<1$.
\eco

\mysection{Basic theorems on regular functions}

Since the Sochocki--Plemelj formula has been built, a large number
of studies on various boundary value problems in the Clifford analysis arise. Obviously,  $\mathcal{E}+S[g]$ is the solution of the jump problem (\ref{22})
for any entire function $\mathcal{E}$.
But ones disregard to prove the fact that they are just all solutions of the jump (\ref{22}) over a long period of time.
This can be ascribed to the simpler Painlev\'{e} problem.

{\bf Painlev\'{e} problem}.
Find a sectionally  regular function $\Phi$,  with $\Gamma$ as its jump surface,  satisfying the boundary value condition
\beq \label{31}
  \Phi^+(x)=\Phi^-(x), \,\,\,\,\,x\in\Gamma.\eeq
\bth[{\rm see \cite{ld}}]\label{T31}
All solutions of Painlev\'{e} problem
 are entire functions.
\eth

This is based on the following theorem.

\bth[Painlev\'{e}  theorem\mbox{\rm \cite{ld}}]\label{T51}
  Let $\Omega\subset\mathcal{R}^{n+1}$ be an open domain and $\Gamma\subset\Omega$ a $C^1$-smooth surface. If $f$ is left $($right$)$ regular in $\Omega\backslash\Gamma$
and $f\in C(\Omega, C(V_n)\big)$, then $f$ is left $($right$)$ regular in $\Omega$.
\eth

Painlev\'{e}  problem is one of the foundations for solving Riemann
boundary value problems, but it has not been seriously treated in Clifford analysis up to now.
Painlev\'{e} theorem problem deals the removable singularity of regular functions on
the whole surface.
  Recently, Luo and Du have given a simple proof of the Painlev\'{e} theorem in \cite{ld},
  which is based upon 3C theorems (Cauchy type integral theorem, Cauchy--Goursat theorem and
  the generalized Cauchy's integral formula). In general,
  the regular functions defined on a region are continuous but not regular to its boundary in discussing its boundary behavior, for example,
$\Phi^+$ and $\Phi^-$ in the jump problem (\ref{22}). So the
generalized Cauchy's theorem and the generalized Cauchy's integral
formula are inevitable in the study of BVPs.

\bde
  Let $\Omega\subset\mathcal{R}^{n+1}$  be a bounded open domain and $\partial\Omega$ a smooth or piecewise smooth surface. $\Omega$ is called a para-ball if $\forall\, f,g\in C\big(\overline{\Omega},C(V_n)\big)$ and $\forall\,\epsilon>0$,  \vspace{1mm}there is an open  domain $\Omega_\epsilon$ with the piecewise smooth surface $\partial\Omega_\epsilon$, such that $\overline{\Omega_\epsilon}\subset\Omega$  and
  \beq
 \left|\displaystyle\int_{\partial\Omega_\epsilon} g(x)\mathrm{d}\sigma f(x)-\int_{\partial\Omega} g(x)\mathrm{d}\sigma f(x)\right|<\epsilon.
 \eeq
  The boundary of a para-ball is called a para-sphere.
\ede

\bth[ Cauchy-Goursat theorem {\rm \cite{ld}}]\label{pseudo-ballCauchytheorem}
 Let $\Omega\subset\mathcal{R}^{n+1}$ be a para-ball.\vspace{1.5mm} If $f$ is  left regular in $\Omega$ and continuous on
 its closure $\overline{\Omega}$, \vspace{1mm}$g$ is right regular in $\Omega$ and continuous on $\overline{\Omega}$, then
 \beq \label{312}
 \displaystyle\int_{\partial\Omega} g(x)\mathrm{d}\sigma f(x)=0.
 \eeq
\eth

\bde
 \label{def42}
 Let $\Omega\subset\mathcal{R}^{n+1}$  be a bounded open domain and $\partial\Omega$ a smooth or piecewise smooth surface
 with the induced orientation by $\Omega$.
$\Omega$ is called a strong  para-ball if for given $w\in\Omega$ and $\forall\,\epsilon>0$,
there are a domain $\Omega_{*}$ with the piecewise smooth surface $\partial\Omega_{*}$ oriented
  by $\Omega_{*}$  and a number $\rho>0$
  such that \beq
 \label{480}
w\in \Omega_{*}^{\circ},\,\,\,\,\,\,\overline{\Omega_{*}}\subseteq\Omega,\,\,\,\,\,\,\partial\Omega_{*}\in U(\partial \Omega,\rho)
 \eeq
  and when $f,g\in C\big(U(\partial\Omega,\rho)\cap\overline{\Omega},C(V_n)\big)$
  \beq\label{481}
 \left|\displaystyle\int_{\partial\Omega_{*}} g(x)\mathrm{d}\sigma f(x)-\int_{\partial\Omega} g(x)\mathrm{d}\sigma f(x)\right|<\epsilon.
 \eeq
  where
  \beq
  U(\partial \Omega,\rho)=\big\{w,\, |w-z|\leq \rho,\, z\in\partial\Omega\big\}
  \eeq
  is a neighbourhood of $\partial\Omega$. The boundary of a strong para-ball is called a strong para-sphere.
\ede

\bth[{\bf Generalized Cauchy's integral formula}]\label{pseudo-ballCauchyformula}
 Let $\Omega\subset\mathcal{R}^{n+1}$ be a strong para-ball.\vspace{0.8mm} If $f$ is  left regular in $\Omega$ and continuous
 on $\overline{\Omega}$, $g$ is right regular in $\Omega$ and continuous on $\overline{\Omega}$, then
 \beq \label{482}
 \displaystyle\frac{1}{\bigvee_{n+1}}\int_{\partial\Omega}E(x-w)\,\mathrm{d}\sigma\, f(x)=
  \left\{ \begin{array}{lll}
    f(w),&when &w\in\Omega,\\[2mm]
    0,&when& w\in\mathcal{R}^{n+1}\backslash\overline{\Omega},
  \end{array}
  \right.
   \eeq
   \beq\label{483+}
 \displaystyle\frac{1}{\bigvee_{n+1}}\int_{\partial\Omega}g(x)\, \mathrm{d}\sigma\,  E(x-w) =
  \left\{ \begin{array}{lll}
    g(w),&when &w\in\Omega,\\[2mm]
    0,&when& w\in\mathcal{R}^{n+1}\backslash\overline{\Omega}.
  \end{array}
  \right.
   \eeq
 \eth

\bre{\rm In fact, Theorem \ref{pseudo-ballCauchytheorem} and Theorem \ref{pseudo-ballCauchyformula} (collectively called 2C theorems) only  schemed out a method to extend weak versions of
  Cauchy's theorem and
 Cauchy's integral  formula (i.e., the classical Cauchy's theorem and
 Cauchy's integral  formula, e.g., see \cite{BDS-RNM-82}) in principle. It
 transfers the key to the difficult question on
  what kind of domains are para-ball or strong para-ball.
  Obviously, the ball $S^{n}$ in $\mathcal{R}^{n+1}$ is a strong para-ball.
  Thus, the Cauchy-Goursat theorem and the generalized Cauchy's
  integral formula hold on $S^n$.  In \cite{zhang13}, Zhang  Zhongxiang
  specifically proved that Cauchy-Goursat theorem holds on $2$-dimensional ball.
  R. Z. Yeh established a kind of  2C theorems  in \cite{Yeh-94}. Those theorems require
  $\Omega$ must be a regular domain in $\mathcal{R}^{n+1},$ even $\Omega$ being cut down a
  small ball is also a regular domain. One thus runs into difficulties
  on classifying regular domains.
 In \cite{Bl00}, to solve
   some quaternionic Riemann boundary value problems,  R. Abreu and J. Bory  presented other kinds
   of generalized
  Cauchy integral formulas involving certain geometric measure theoretic
  concepts. Those generalized
  Cauchy integral formulas cannot be deduced directly from the Gauss-Green theorem in geometric
measure theory unless appending the condition that each component of
$D[f]$ is absolutely integrable on $\Omega$.
  }\ere

 In the following we shall illustrate that the domains bounded by smooth surfaces and a special type of piecewise smooth surfaces are strong para-balls.

\bex\label{example41}
{\rm
Let $\Omega\subset\mathcal{R}^{n+1}$  be a bounded open domain and its boundary $\partial\Omega$ a $C$-smooth surface with the induced orientation by $\Omega$. Then $\Omega$ is a
strong  para-ball.
}
\eex

\bde Fix $\ell\in \{0,1,\cdots,n\}$.
Let $\varphi^{+}$ and $\varphi^{-}$ be  elementary $C^{1}$-smooth surfaces in $\mathcal{R}^{n+1}$
with the explicit representations
\beq\label{e32}
x_\ell=\varphi^{\pm}\big(u^{\ell}\big),\,\,u^{\ell}=\big(x_0,\cdots,x_{\ell-1},x_{\ell+1},\cdots,x_{n}\big)\in U_{\ell}
\eeq
where the parametric
 domain $U_\ell$  is an open domain in  $\mathcal{R}^{n}$ covering the cube
 \beq
 I^{\ell}=\big\{u^{\ell},
 \,\,a^{-}_k\leq x_k\leq a^{+}_k,\, \,k=0,\cdots\!,\ell-1,\ell+1,\cdots,n\big\}.
 \eeq
If $\varphi^{-}\big(u^{\ell}\big)<\varphi^{+}\big(u^{\ell}\big)$  for  $u^{\ell}\in U_\ell$, we call
\beq
\label{494}
\Pi^{\ell}\!=\!\Big\{x,\,\,x\!=\!\big(x_0,x_1,\cdots,x_n\big)\in\mathcal{R}^{n+1},\,\,\varphi^{-}\big(u^{\ell}\big)\!<\!x_{\ell}\!<\!\varphi^{+}\big(u^{\ell}\big),\,\,a^{-}_j\!<\!x_j\!<\!a^{+}_j,\,\, j\!\not=\!\ell\Big\}
 \eeq
the  $x_{\ell}$-type cylindrical body with crooked tips in $\mathcal{R}^{n+1}$.
 \ede

\bex\label{e23}
Let $\Pi^{\ell}$ be a $x_{\ell}$-type cylindrical body  with crooked tips in $\mathcal{R}^{n+1}$.
If its boundary $\partial\Pi^{\ell}$ is oriented by the exterior normal, then it is a strong para-ball.
\eex

The authors in \cite{ld} give the two examples in detail. Now we show some useful cases of $2$C theorem.

\bth[Cauchy-Goursat theorem on  smooth surface {\rm \cite{ld}}] \label{smoothsurface1}
  Let $\Omega\!\subset\!\mathcal{R}^{n+1}$ be a compact domain and $\partial\Omega$ a $C^{1}$-smooth surface.
  If $f$ is left regular and $g$ is right regular in $\Omega^\circ$ and $f,g\in C(\Omega, C(V_n))$, then
  \beq
  \int_{\partial\Omega} g(x)\,\mathrm{d}\sigma f(x)=0.
  \eeq
\eth

\bth[Generalized Cauchy integral formula on  smooth surface {\rm \cite{ld}}] \label{smoothsurface1 2}
  \vspace{0.5mm} Let $\Omega\subset\mathcal{R}^{n+1}$ be a bounded open domain and $\partial\Omega$ a $C^{1}$-smooth surface.
 If $f$ is  left regular in $\Omega$ and continuous on $\overline{\Omega}$, $g$ is right regular in $\Omega$ and continuous on $\overline{\Omega}$, then
\beq
 \displaystyle\frac{1}{\bigvee_{n+1}}\int_{\partial\Omega}E(x-w)\,\mathrm{d}\sigma\, f(x)=
  \left\{ \begin{array}{lll}
    f(w),&when &w\in\Omega,\\[2.5mm]
    0,&when& w\in\mathcal{R}^{n+1}\backslash\overline{\Omega},
  \end{array}
  \right.
 \eeq
\beq
   \displaystyle\frac{1}{\bigvee_{n+1}}\int_{\partial\Omega}g(x)\, \mathrm{d}\sigma\, E(x-w) =
  \left\{ \begin{array}{lll}
    g(w),&when &w\in\Omega,\\[2.52mm]
    0,&when& w\in\mathcal{R}^{n+1}\backslash\overline{\Omega}.
  \end{array}
  \right.
 \eeq
\eth

\bth[2C theorem on  cylindrical body  {\rm \cite{ld}}] \label{cylindrical}
  Let $\Omega$ be an $x_{\ell}$-type cylindrical body with crooked tips in $\mathcal{R}^{n+1}$ and the orientation of its boundary $\partial \Omega$
   induced  by the
  exterior normal.\vspace{1mm}
    If $f$ is left regular and $g$ is right regular in $\Omega$ and $f,g\in C\big(\overline{\Omega}, C(V_n)\big)$, then
 \beq
 \label{314}
 \int_{\partial\Omega} g(x)\,\mathrm{d}\sigma f(x)=0,
 \eeq
 \beq
 \displaystyle\frac{1}{\bigvee_{n+1}}\int_{\partial\Omega}E(x-w)\,\mathrm{d}\sigma\, f(x)=
      f(w),\,\,w\in\Omega,
\eeq
\beq \label{315}
 \displaystyle\frac{1}{\bigvee_{n+1}}\int_{\partial\Omega}g(x)\, \mathrm{d}\sigma\,  E(x-w) =
   g(w),\,\,\,w\in\Omega.
\eeq
  \eth

  From the  Sochocki--Plemelj formula (\ref{231}), (\ref{314}) and (\ref{315}) we know that the following lemma holds.
\ble\label{317}
If $f\in H(\Gamma)$, then it is the positive boundary value  of a left regular function $F$
 on $\Omega^{+}$, if and only if one of the conditions
\beq
\frac{1}{2}f(t)=\Big(C[f]\Big)(t)=\displaystyle\frac{1}{\bigvee_{n+1}}\int\limits_{\Gamma}
 \frac{\overline{x-t}}{|x-t|^{n+1}}\,\mathrm{d}\sigma f(x), \hspace{4mm}t\in\Gamma.
\eeq
and
\beq
\Big(C[f]\Big)(w)=\displaystyle\frac{1}{\bigvee_{n+1}}\int\limits_{\Gamma}
 \frac{\overline{x-w}}{|x-w|^{n+1}}\,\mathrm{d}\sigma f(w)=0, \hspace{4mm}w\in\Omega^{-}
\eeq
is satisfied.
\ele

  Now we give another version of the Cauchy integral formula.

  \bth[Version of Cauchy's integrals theorem in span form]\label{T28}
If $f$ is left regular on $\Omega^+$ and H\"{o}lder continuous to  $\Gamma$, then
\beq\label{318}
\Big(C[f]\Big)(w)=\displaystyle\frac{1}{\bigvee_{n+1}}\int\limits_{\Gamma}
 \frac{\overline{x-w}}{|x-w|^{n+1}}\,\mathrm{d}\sigma f(x)=\mbox{\rm span}(\Gamma,w)\,f(w),
 \hspace{4mm}w\in\mbox{$\cal R$}^{n+1}
\eeq
where
\beq\label{319}
\mbox{\rm span}(\Gamma,w)=\left\{
\begin{array}{ll}
1, &w\in\Omega^+,\\[2mm]
\displaystyle\frac{1}{2}, &w\in\Gamma,\\[4mm]
0, &w\in\Omega^{-}
\end{array}
\right.
\eeq
is the \mbox{\rm span} at $w$ with respect to $\Gamma$.
\eth

\bre{\rm
The span about $\Omega$ at $w$ may be visually understood as follows. For $w\in \Omega^+$, since its entire neighborhood lies in $\Omega^+$, so
$\mbox{\rm span}=1$. For $w\in \Omega^{-}$, since its entire neighborhood lies in  $\Omega^{-}$,  so
$\mbox{\rm span}=0$.
For $w\in\Gamma$,  its neighborhood partly lies both $\Omega^+$ and $\Omega^{-}$, moreover $\Gamma$ is smooth,  so $\mbox{\rm span}=\frac{1}{2}$.
}
\ere

To discuss behavior of a regular function, now we introduce its
order at the infinity.

\bde
Let $\Phi$ be  regular on some neighborhood of $\infty.$ If there exists an integer $m$ such that
  \begin{equation}
    0<\limsup\limits_{|x|\rightarrow +\infty}|x|^{-m}|\Phi(x)|<+\infty
  \end{equation}
  then $m$ is called the order of $\Phi$ at $\infty$, that is
 denoted by $\mbox{\rm Ord\,}(\Phi,\infty)=m.$ When $f\equiv 0,$ we have that $\mbox{\rm Ord} (f,\infty)=-\infty.$
\ede

In order to deal with the order issue at the infinity  of regular functions,
we need to review Taylor expansion of entire functions and expansion of the Cauchy type integral.
For the Taylor series and the Laurent series theory in Clifford analysis, we refer  to
 \cite{BDS-RNM-82, Delanghe-70, Delanghe-72, Yeh-91, Yeh-94 }.
Let $Z=(z_1,\cdots,z_n)$ and $\alpha=[\alpha_1,\alpha_2,\cdots,\alpha_n],$ where
$z_{j}$'s are
hypercomplex variables given by
\beq
z_j=z_j(x)=x_je_0-x_0e_j,\,\,\, \,\,x\in\mathcal{R}^{n+1},\,\,\,\,\,\,\,\, ( j=1,\cdots,n).
\eeq
Note that the functions $z_j$ are biregular \cite{BDS-RNM-82,Delanghe-70,Delanghe-72,Yeh-94,Yeh-91, ts},
and $\alpha_j$'s are nonnegative integers. The symmetry power $Z^\alpha$ is a
biregular function in $\mathcal{R}^{n+1},$ defined as the sum of all possible $z_i$ products, each of which contains $z_i$ factor exactly $\alpha_i$ times.
For example, for $n=2$
\beq
(z_1,z_2)^{[0,0]}=1,\,\,\,
(z_1,z_2)^{[1,1]}=z_1z_2+z_2z_1,\,\,\,(z_1,z_2)^{[2,0]}=z_1^2.
\eeq

Let $N_0$ stand for the set
$\{0, 1,2,\cdots, \}$. Write
\beq
\big|\alpha\big|=\sum\limits_{j=1}^{n}\alpha_j
\,\,\,\mbox{for}\,\,\,\alpha=[\alpha_1,\alpha_2,\cdots,\alpha_n]\in N_0^{n},
\eeq
\beq
\partial^{\alpha}=\displaystyle\frac{\partial^{|\alpha|}}
{\partial x_1^{\alpha_1}\cdots\partial{x_n^{\alpha_n}}}.
\eeq

\bth[Taylor's expansion {\rm \cite{BDS-RNM-82,Yeh-94,Yeh-91,ts}}]\label{Laurentseries}
  Let $f$ be left entire in $\mathcal{R}^{n+1}$. Then  $f$ can be expanded into a Taylor series
  \beq \label{320}
    f(x)=\sum\limits_{k=0}^\infty P_k[f](x),\,\,\,\,\, x\in\mathcal{R}^{n+1}
  \eeq
  where \beq
    P_k[f](x)=\displaystyle\frac{1}{k!}\sum_{|\alpha|=k}Z^\alpha(x)\big[\partial^{\alpha}f\big](0).
   \eeq
   Similarly, if $f$ is right entire, then
 \beq\label{321+}
    P_k^{r}[f](x)=\displaystyle\frac{1}{k!}\sum_{|\alpha|=k}\big[\partial^{\alpha}f\big](0)Z^\alpha(x).
   \eeq
 \eth

\bre{\rm Since
\beq \label{321}
\big[\partial^{\alpha}f\big](0)=\frac{(-1)^{|\alpha|}}{\bigvee_{n+1}}
\int_{\partial B(0,R)}\big[\partial^{\alpha}E\big](x)\,\mathrm{d}\sigma\, f(x),\,\,\,R\in (0,+\infty),
  \eeq
we have
\beq \label{322}
\Big|P_k[f](x)\Big|\leq M\, C_{k+n-1}^{k+1}\, \Big(\!1+k^2\Big)\,
\displaystyle\frac{|x|^k}{R^k}\,\,\sup\Big\{\big|f(y)\big|,\,\,|y|=R\Big\},\,\,\,\,\,\,|x|<R,
    \eeq
where the constant $M>0$ only depends upon the dimension $(n+1)$.
   }
\ere

We immediately obtain the following important result through (\ref{320})-(\ref{322}).

\bth[Liouville type theorem]\label{Liouvilletheorem}
If $f$ is  entire with $\mbox{\rm Ord}(f,\infty)\leq m$, then it is a
hyperpolynomial of  degree $m$. For example, if $f$ is left entire, then
\beq \label{325}
f(x)=\left\{
\begin{array}{lll}
\displaystyle\sum_{|\alpha|=0}^{m}\displaystyle\frac{1}{|\alpha|!}\, Z^{\alpha}(x)\,c_{\alpha}, &\mbox{when}\,&m\geq 0,\\[6mm]
0,&\mbox{when}\,&m<0,
\end{array}\right.
\eeq
 where $c_{\alpha}$'s are hypercomplex constants.
\eth

\bde
We call
\beq
f(x)=\sum_{|\alpha|=0}^{m} Z^{\alpha}(x)\,c_{\alpha}\,\, \mbox{with}\,\,
\sum_{|\alpha|=m}|c_{\alpha}|\not=0
\eeq
 a hypercomplex symmetric  polynomial of degree $m$, denoted by $\mbox{\rm deg}(f)=m$.
\ede

\bex[\hspace{-0.2mm}{\rm \cite{ld}}]\label{e33}
A hypercomplex polynomial $P_m$ is of degree $m$ if and only if $\mbox{\rm ord}(P_m,\infty)=m$.
\eex

\bth[Laurent's expansion {\rm \cite{BDS-RNM-82}}]\label{Laurentseries+} The Cauchy type integral
$S[g]$ has the Laurent series expansion near the infinity
\beq \label{641+}
\Big(S[g]\Big)(w)=-\sum\limits_{k=0}^{\infty} Q_k(w),\,\,\,\,\,\,w>\rho=
\max\big\{|x|,\, x\in\Gamma\big\},
\eeq
where
\beq
Q_k(w)=\sum\limits_{|\alpha|=k}\big[\partial^{\alpha}E\big](w)\frac{1}{|\alpha|!}\int_{\Gamma}Z^\alpha(x)\mathrm{d}\sigma g(x).
\eeq
moreover,
\beq \label{638+}
\Big|Q_{k}(w)\Big|\leq M \,C_{k+n-1}^{k+1}\, \Big(1+k^2\Big)\, \frac{\rho^k}{\big|w|^{n+k}}\displaystyle
\int_{\Gamma}\big|g(x)|\mbox{\rm d}s,\,\,\,|w|>\rho,
 \eeq
  where the constant $M>0$ only depends upon the dimension $n$. Similarly, The Cauchy type integral
$[g]S$ has the Laurent series expansion near the infinity
\beq \label{641}
\Big([g]S\Big)(w)=-\sum\limits_{k=0}^{\infty} Q_k^{r}(w),\,\,\,\,\,\,w>\rho=
\max\big\{|x|,\, x\in\Gamma\big\},
\eeq
where
\beq
Q_k^{r}(w)=\sum\limits_{|\alpha|=k}\big[\partial^{\alpha}E\big](w)\frac{1}{|\alpha|!}\int_{\Gamma} g(x)\,\mathrm{d}\sigma \,Z^\alpha(x).
\eeq
moreover,
\beq \label{638}
\Big|Q_{k}^{r}(w)\Big|\leq M \,C_{k+n-1}^{k+1}\, \Big(1+k^2\Big)\, \frac{\rho^k}{\big|w|^{n+k}}\displaystyle
\int_{\Gamma}\big|g(x)|\mbox{\rm d}s,\,\,\,|w|>\rho,
 \eeq
  where the constant $M>0$ only depends upon the dimension $n$.
\eth

  \ble[Order of  Laurent hypercomplex polynomial {\rm \cite{ld}}]
\label{l31}
Let
\beq
Q_m^{l}(x)=\sum_{|\alpha|=m}\big[\partial^{\alpha} E\big](x)\lambda_{\alpha},\,\,\,\,Q_{m}^{r}(x)=\sum_{|\alpha|=m}\lambda_{\alpha}\big[\partial^{\alpha} E\big](x)
\eeq
be the hypercomplex Laurent  polynomials
where  $\lambda_{\alpha}$'s are some hypercomplex constants. Then we have that $\mbox{\rm Ord}\big(Q_m^{l},\infty\big)\!=\!-n\!-\!m$ or
$\mbox{\rm Ord}\big(Q_m^{r},\infty\big)\!=\!-n\!-\!m$ if and only if
\beq \label{640}
\sum\limits_{|\alpha|=m}|\lambda_{\alpha}|\not=0.
\eeq
\ele

\bex[Order of Cauchy type integral {\rm \cite{ld}}]\label{e34}
Let
\beq \label{337}
N^{l}=\min\left\{|\alpha|,\,\displaystyle\int_{\Gamma}Z^\alpha(x)\mathrm{d}\sigma g(x)\not=0,\, \alpha\in N_0^{n}\right\},
\eeq
\beq
N^r=\min\left\{|\alpha|,\,\displaystyle\int_{\Gamma} g(x)\,\mathrm{d}\sigma \,Z^\alpha(x)\not=0,\, \alpha\in N_0^{n}\right\},
 \eeq
 provided
the above sets are not empty,
then
\beq
\mbox{\rm Ord}(S[g],\infty)\!=\!-n\!-\!N^{l},\,\,\,\,\mbox{\rm Ord}([g]S,\infty)\!=\!-n\!-\!N^r.
\eeq
\eex

\mysection{Riemann boundary value problems}
The jump problem (\ref{22}) with the restrictions on the growth condition  at the infinity for sectionally regular functions give rise to the following problem.
\vspace{1mm}

{\bf $R_m$ jump problem.}
Find a sectionally  regular function $\Phi$,  with $\Gamma$ as its jump surface,  satisfying the boundary value condition and
the growth condition at the infinity

\beq \label{41}
 \left\{
   \begin{array}{l}
 \Phi^+(x)=\Phi^-(x)+g(x),\,\,\,\, x\in\Gamma,\\[2mm]
 \mbox{\rm Ord\,}(\Phi, \infty)\leq m,
 \end{array}
 \right.
 \eeq
 where $g\in H(\Gamma)$. \vspace{1mm}

The simplest $R_m$ problem is the case of ({\ref{41}) when $g=0$, which is called Liouville problem.\vspace{1mm}

{\bf Liouville problem.}
Find a sectionally  regular function $\Phi$,  with $\Gamma$ as its jump surface,  satisfying the boundary value condition and
the growth condition at the infinity
\beq \label{42}
 \left\{
   \begin{array}{l}
 \Phi^+(x)=\Phi^-(x),\,\,\,\, x\in\Gamma,\\[2mm]
 \mbox{\rm Ord\,}(\Phi, \infty)\leq m.
 \end{array}
 \right.
 \eeq

Obviously, any hypercomplex symmetric polynomials of degree not exceeding $m$ is
the solution of Liouville problem (\ref{42}). On the other hand, any solution of Liouville problem (\ref{42}) is a symmetric polynomial of degree not exceeding $m$
by Theorem \ref{T31} and  Example \ref{e33}.

\ble\label{l41}
The solutions of the Liouville problem $(\ref{42})$ are all hypercomplex symmetric polynomials of degree not exceeding $m$.\vspace{-2mm}
\ele

 Below we firstly solve  jump $R_{-n}$ problem, i.e., (\ref{41}) with $m=-n,$ that is a fundamental case of $R_m$ problem,
 then the general $R_{m}$ jump problem.\vspace{1mm}

\mbox{\bf Jump $R_{-n}$ problem.} Find a  sectionally regular function $\Phi$, with $\Gamma$ as its jump surface, such that
\beq\label{jp-n} \left\{
\begin{array}{l}
\Phi^{+}(x)-\Phi^{-}(x)=g(x),\,\, x\in\Gamma,\\[2mm]
\mbox{\rm Ord\,}(\Phi, \infty)\leq -n,\end{array}
\right.
\eeq
 where $g\in H(\Gamma)$. \vspace{2mm}

  By the Sochocki-Plemelj formulae and Example \ref{e34}, we know that $S[g]$ $([g]S)$ is just  the solution of the jump $R_{-n}$ problem.
Let $\Delta=S[g]-\Psi$, then  $\Psi$ is the solution of the $R_{-n}$ problem if and only if
$\Delta$ is the solution of the
  Liouville problem (\ref{42}) with $m=-n.$
 By Lemma \ref{l41}, we get $\Delta=0,$ that is to say that $S[g]$ is the unique solution of the jump $R_{-n}$ problem. Then  we have the following Lemma.
  \ble \label{l42}
  The $R_{-n}$  jump problem $(\ref{jp-n})$ has the unique  solution $S[g]$ $([g]S)$.\vspace{-2mm}
  \ele

Now we discuss the $R_m$ problem $(\ref{41})$, which is divided into three cases. Below we only find the left regular
solution.

  \mbox{\bf Case $1$}: $m\geq 0$. Firstly, $S[g]$ is a solution by Lemma \ref{l42}. Secondly, letting $\Delta=S[g]-\Psi$, then $\Psi$ is the solution of the
  $R_{m}$ problem  if and only if
$\Delta$ is the solution of the
  Liouville problem (\ref{42}).
 By quoting Lemma \ref{l41}, we get $\Delta=P_m,$ that is arbitrary hypercomplex polynomial of degree not exceeding $m$ denoted by $P_m\in \prod_m^{H}.$ Thus,
 the general solution of the $R_m$ problem $(\ref{41})$ is
  \beq \label{44}
 \Phi(w)=\Big(S[g]\Big)(w)+ \sum_{|\alpha|=0}^{m}Z^{\alpha}\,c_{\alpha}
 \eeq
 where $c_{\alpha}$'s are $C^{m}_{n+m}$ arbitrary hypercomplex constants.

  {\bf Case $2$}: $-n\leq m<0$.  $S[g]$ is a solution by Lemma \ref{l42}.
      Similarly, $\Psi$ is the solution of the jump $R_{m}$ problem if and only if
$\Delta$ is the solution of the
  Liouville problem (\ref{42}). Hence,
 quoting Lemma \ref{l41}, we have $\Delta=0$. Thus,
  the $R_m$ problem $(\ref{41})$ has unique solution $S[g]$.

{\bf Case $3$}: $m<-n$.  By Lemma \ref{l42}, if the jump problem is solvable then its solution must be the solution of $R_{-n}$ problem $S[g]$.
 Conversely,  if $S[g]$ is the solution of the  jump problem $R_{-n}$, then,  by Example \ref{e34}
  it  must satisfy the conditions of solvability
   \beq \label{45}
         \int_{\Gamma}Z^{\alpha}(x)\mbox{\rm d}\sigma g(x)=0,\,\,\,|\alpha|=0,1,\cdots\!, -(n-1+m).
 \eeq

 Summarizing the above discussions, we have the following result.

\bth\label{jump}
 For the jump problem $R_m$ $(\ref{41})$, the general solution is stated  in three cases according to values of $m$.

   $(1)$  When $m\geq 0$, the jump problem is solvable, its general solution is given by $(\ref{44})$ with
   $C^{m}_{n+m}$ free hypercomplex constants.

 $(2)$  When $-n\leq m<0$, the jump problem has a unique solution $S[g]$.

$(3)$  When  $m<-n$, the jump problem has a unique solution $S[g]$
if and only if the $C^{n}_{-m-1}$ conditions of solvability in  $(\ref{45})$
are satisfied.
  \eth

  \bre{\rm  The number of $\alpha$ with $|\alpha|=k$ is $C^{n-1}_{k+n-1}=C^{k}_{k+n-1}.$ \vspace{1.5mm}Hence the number of $c_{\alpha}$'s in
  (\ref{44}) is just $C_{n+m}^{m}$. Similarly,
 (\ref{45}) includes $C^{n}_{-m-1}$ conditions of solvability.
  }\ere

\bre{\rm
 There is no  regular function $F$ near the infinity such that $\mbox{\rm Ord}(g,\infty)=m$ with $-n<m<0$.
In fact, if  $g(w)$ is regular on $|w|> R$, taking
\beq
\Gamma=\{w, |w|=2R\}
\eeq
then $g\in H(\Gamma)$
 and
\beq
\Phi(w)=\left\{
\begin{array}{ll}
0,&|w|<2R,\\[2mm]
g(w), &|w|>2R
\end{array}
\right.
\eeq
is the solution of the jump $R_{-n}$ problem  $(\ref{jp-n})$, so $\Phi=S[g]$ by Theorem \ref{jump}, which results in $\mbox{\rm Ord}(g,\infty)\leq-n$.
 This fact is very interesting, since there is analytic function near the infinity with arbitrary order in the classical complex analysis. }
\ere

{\bf  $R_m$ problem with constant gap $G$}.
Find a sectionally  left regular function $\Phi$,  with the jump surface $\Gamma$,  satisfying the boundary value condition and
the growth condition at the infinity
\beq \label{47}
 \left\{
   \begin{array}{l}\Phi^+(x)=\Phi^-(x)G+g(x),\,\,\,\, x\in\Gamma,\\[2mm]
 \mbox{\rm Ord\,}(\Phi, \infty)\leq m,
 \end{array}
 \right.
 \eeq
where $g\in H(\Gamma)$ and $G$ is an invertible hypercomplex constant.\vspace{2mm}

Let $\Psi=\Phi X^{-1}$, the $R_m$ problem $(\ref{47})$ becomes a $R_m$ jump problem for $\Psi$. So, we obtain immediately the following result.

\ble \label{lemma68}
The general solution of the
$R_m$ problem $(\ref{47})$ is
\beq \label{48}
\Phi(w)
=\left[\Big(S\big[g\big]\Big)(w)+
 P_m(w)\right]\!X(w)
 ,\hspace{3mm}w\not\in \Gamma,
\eeq
where
\beq\label{49}
X(w)=\left\{
\begin{array}{ll}
1, &w\in\Omega^+,\\[2mm]
G^{-1}, &w\in\Omega^{-},
\end{array}
\right.
\eeq
and $P_m$ is an arbitrary hypercomplex symmetric polynomial of degree not  exceeding $m$.
When  $m<-n$, then the jump problem has a unique solution $(\ref{48})$ with $P_m=0$
if and only if the following $C^{n}_{-m-1}$ conditions of solvability
         \beq \label{410}
         \int_{\Gamma}Z^{\alpha}(x)\mbox{\rm d}\sigma g(x)=0,\,\,\,|\alpha|=0,1,\cdots\!, -(n-1+m)
 \eeq
are fulfilled.
\ele

\bre{\rm In \cite{Duzh-01}, since the order at the infinity is not well defined, Zhang and Du only discuss the $R$ problem, i.e., the growth condition at the infinity of (\ref{47}) is ignored. Still, the result there must be based on the Painlev\'{e} problem \ref{23},
which is not discussed in earlier works. In \cite{du+}, we reproved the Painlev\'{e} problem by the version of Cauchy's integral theorem in span form (\ref{318}),
 but the latter is based on the generalized 2C theorems, that is still
not proved there.
}
\ere

\rm {\bf General $R_m$  problem.}
Find a sectionally  left regular function $\Phi$,  with $\Gamma$ as its jump surface,  satisfying the boundary value condition and
the growth condition at the infinity

\beq \label{411}
 \left\{
   \begin{array}{l}
 \Phi^+(x)=\Phi^-(x)G(x)+g(x),\,\,\,\, x\in\Gamma,\\[3mm]
 \mbox{\rm Ord\,}(\Phi, \infty)\leq m,
 \end{array}
 \right.
 \eeq
 where $G,g\in H(\Gamma)$.

\bre{ The general $R_m$  problem $(\ref{411})$ is an open problem up to now. We believe that there will arise new tools to solve this problem.
}
\ere

\mysection{Singular integral equations}

Singular integral equations of higher dimensions are frequently encountered in physical and engineering applications. There is
very close link between boundary value problems and singular integral equations. We consider the singular
integral equation
\beq \label{51}
\varphi(t)a(t)+\displaystyle\frac{2}{\bigvee_{n+1}}\int\limits_{\Gamma}
 \frac{\overline{x-t}}{|x-t|^{n+1}}\,\mathrm{d}\sigma \varphi(x)k(x,t)=f(t),\hspace{4mm}t\in\Gamma,
\eeq
where the input function $k\in H(\Gamma\times\Gamma)$, $a, f\in H(\Gamma)$ and the output function
$\varphi\in H(\Gamma)$. Since the operator
\beq \label{52}
\Big(\mathcal{K}[\varphi]\Big)(t)=\displaystyle\frac{2}{\bigvee_{n+1}}\int\limits_{\Gamma}
 \frac{\overline{x-t}}{|x-t|^{n+1}}\,\mathrm{d}\sigma \varphi(x)\big[k(x,t)-k(t,t)\big]
 \eeq
 is compact, we will consider the characteristic equation
 \beq \label{53}
 \varphi(t)a(t)+\displaystyle\frac{2}{\bigvee_{n+1}}\int\limits_{\Gamma}
 \frac{\overline{x-t}}{|x-t|^{n+1}}\,\mathrm{d}\sigma \varphi(x)b(t)=f(t),\hspace{4mm}t\in\Gamma.
\eeq
Let
\beq \label{54}
\Phi(z)=\displaystyle\frac{1}{\bigvee_{n+1}}\int\limits_{\Gamma}
 \frac{\overline{x-z}}{|x-z|^{n+1}}\,\mathrm{d}\sigma\varphi(x),\hspace{4mm}z\not\in\Gamma.
\eeq
Then,  from (\ref{54})  we get the Riemann boundary problem
\beq \label{55+}
\left\{
\begin{array}{ll}
\Phi^{+}(t)\big[a(t)+b(t)\big]=\Phi^{-}(t)\big[a(t)-b(t)\big]+f(t), \hspace{4mm}t\in\Gamma,\\[3mm]
\Phi(\infty)=0.
\end{array}
\right.
\eeq

\ble \label{l51}
If $\varphi$ is the solution of the singular integral equation $(\ref{53}),$ then $\Phi$ given in $(\ref{54})$
is the solution of the Riemann boundary problem $(\ref{55+})$. Inversely, if $\Phi$ is
the solution of Riemann boundary problem $(\ref{55+}),$ then $\varphi(t)=\Phi^{+}(t)-\Phi^{-}(t)$
is the solution of the singular integral equation $(\ref{53})$.
\ele

Due to obstacles not being overcome in the general $R_m$  problem $(\ref{55+}),$
the general case of (\ref{53}) is also an open problem, a special case
 of which is
solved. Assume that both $a(t)+b(t)$ and $a(t)-b(t)$ are inversable and
\beq\label{55}
\displaystyle\frac{a(t)-b(t)}{a(t)+b(t)}\!|=G
\eeq
 is
a hypercomplex constant, then the singular integral equation ({\ref{53}) is
reduced to an $R_{-1}$ problem as follows.
\beq \label{56}
\left\{
\begin{array}{ll}
\Phi^+(t)=\Phi^{-}(t)\displaystyle\frac{a(t)-b(t)}{a(t)+b(t)}\!|+
\displaystyle\frac{f(t)}{a(t)+b(t)}\!|\, , \hspace{4mm}t\in\Gamma,\\[4mm]
\Phi(\infty)=0.
\end{array}
\right.\eeq

By Lemma \ref{l51} and Theorem \ref{jump}, we have
\bth
Singular integral equation $(\ref{53})$ has an unique solution
\beq \label{58}
\phi(t)\!=\!\displaystyle\frac{f(t)}{2}\!\left(\frac{1}{a(t)\!+\!b(t)}+\frac{1}{a(t)\!-\!b(t)}\right)
\!-\!\displaystyle\frac{2}{\bigvee_n}\int\limits_{\Gamma}
 \frac{\overline{x-t}}{|x-t|^n}\,\mathrm{d}\sigma\frac{f(x)}{a(x)-b(x)}\!|\,
 \frac{b(x)}{a(x)+b(x)}\!|,\,\,\,\,
t\in\Gamma.
\eeq
\eth

While $a=0$, $b=1$, we get the inversion formula for Cauchy principal value integral.
\bth
[Inversion formula for Cauchy principal value integral\mbox{\rm \cite{Duzh-01,ld}}] If $f\!\in\!H(\Gamma)$, then
\beq \label{59}
f(t)=\displaystyle\frac{1}{\bigvee_{n+1}}\int\limits_{\Gamma}
 \frac{\overline{x-t}}{|x-t|^{n+1}}\,\mathrm{d}\sigma\varphi(x)\Longleftrightarrow
\varphi(t)=\displaystyle\frac{1}{\bigvee_{n+1}}\int\limits_{\Gamma}
 \frac{\overline{x-t}}{|x-t|^{n+1}}\,\mathrm{d}\sigma f(x).
\eeq
 \eth

\bre{\rm The inversion formula for Cauchy principal value integral may also be directly obtained from Theorem \ref{T28}.
It shows the singular integral operator
\beq \label{510}
\Big(S[f]\Big)(t)=\displaystyle\frac{1}{\bigvee_{n+1}}\int\limits_{\Gamma}
 \frac{\overline{x-t}}{|x-t|^{n+1}}\,\mathrm{d}\sigma f(x), \hspace{4mm}t\in\Gamma
\eeq
is idempotent, i.e., $S^2=\mbox{\bf I}$ where $\mbox{\bf I}$ is the identity operator. Moreover,
from
\beq \label{511}
f(t)=\displaystyle\frac{1}{\bigvee_{n+1}}\int\limits_{\Gamma}
 \frac{\overline{x-t}}{|x-t|^{n+1}}\,\mathrm{d}\sigma_x \displaystyle\frac{1}{\bigvee_{n+1}}\int\limits_{\Gamma}
 \frac{\overline{\tau-x}}{|\tau-x|^{n+1}}\,\mathrm{d}\sigma_{\tau} f(\tau),\hspace{5mm}t\in\Gamma,
\eeq
and \cite{ld}
 \beq\label{512}
 \displaystyle\int_{\Gamma}E(x-t) \mathrm{d}\sigma\,  E(\tau-x) =0,\,\,\,\,\,\tau,\, t\in\Gamma,\,\,\,\tau\not=t,
   \eeq
we see that the inversion formula for Cauchy principal value integral is just  the
Poincar\'{e}--Bertrand formula of $f$, i.e.,
\beq
\label{513}
\begin{array}{ll}
&\displaystyle\int_{\Gamma}
 \frac{\overline{x-t}}{|x-t|^{n+1}}\,\mathrm{d}\sigma_{x}\!\left[\,
\displaystyle\int_{\Gamma}
 \frac{\overline{\tau-x}}{|\tau-x|^{n+1}}\,\mathrm{d}\sigma_{\tau} f(\tau)\right]\\[6mm]
 =\!&\!\!\left[\displaystyle\frac{\bigvee_{n+1}}{2}\right]^2 f(t)+
 \displaystyle\int_{\Gamma}\left[\,\displaystyle\int_{\Gamma}
 \frac{\overline{x-t}}{|x-t|^{n+1}}\,\mathrm{d}\sigma_{x}
 \frac{\overline{\tau-x}}{|\tau-x|^{n+1}}\right]\!
\mathrm{d}\sigma_{\tau}\, f(\tau),\,\,\,t\in\Gamma.
 \end{array}
 \eeq
The general Poincar\'{e}--Bertrand formula

  \beq\label{514}
\begin{array}{ll}
&\displaystyle\int_{\Gamma}
 \frac{\overline{x-t}}{|x-t|^{n+1}}\,\mathrm{d}\sigma_{x}\!\left[\,
\displaystyle\int_{\Gamma}
 \frac{\overline{\tau-x}}{|\tau-x|^{n+1}}\,\mathrm{d}\sigma_{\tau} k(\tau,x)\right]-\left[\displaystyle\frac{\bigvee_{n+1}}{2}\right]^2 k(t,t)\\[8mm]
 =&\!\!\sum\limits_{k=0}^{n}(-1)^{k}
 \displaystyle\int_{\Gamma}\left[\,\displaystyle\int_{\Gamma}
 \frac{\overline{x-t}}{|x-t|^{n+1}}\,\mathrm{d}\sigma_{x}
 \frac{\overline{\tau-x}}{|\tau-x|^{n+1}}e_k
 k(\tau,x)\right]\mbox{\rm d}
\widehat{\tau}_{k},\,\,\,\,\,t\in\Gamma
 \end{array}
  \eeq
  is  an open problem,
which will be an important tool for  discussing singular integral equations. We hope that the following conjecture is true, but
 its proof is not straightforward.}
 \ere

\mbox{\bf Conjecture [Poincar\'{e}--Bertrand formula].} {\it If $k\in H(\Gamma\times\Gamma)$, then $(\ref{514})$ holds.}\vspace{2mm}

Other studies for singular integral equations can be found in the mentioned references.

\mysection{Dirichlet boundary value problems}

Throughout this section, we assume that $\Omega=\Omega^{+}$ is a bounded open domain in $\mathcal{R}^{n+1}$ with a $C^1$-smooth boundary $\Gamma$,
 $\Gamma$ is oriented by exterior normal.
  Regular (entire) function means left regular (entire) function. All results for the right regular (entire) function can be similarly obtained. If $G$ in the boundary condition of (\ref{47}) equals to $0$ and the regular function $\Phi$ that is found is in $\Omega,$ then the boundary value problem will become to the Dirichlet boundary value problem. There are several cases based on different continuity of given function.
\vspace{2mm}

{\bf Dirichlet  problem with $H$ continuous boundary  function.}
Find a left regular function $\Phi$ on $\Omega$,  satisfying the boundary value condition\vspace{-1mm}
\beq \label{61}
  \Phi^+(x)=g(x),\,\,\,\, x\in\Gamma,\eeq
where $g\in H(\Gamma)$ is given. \vspace{1mm}

This problem will be simply solved. Firstly, if $\Phi$ is the solution of the problem (\ref{61}), by the generalized Cauchy's integral formula in Theorem \ref{smoothsurface1 2}  we have

\beq\label{62}
\Phi(w)=\big(S[g]\big)(w)=\displaystyle\frac{1}{\bigvee_{n+1}}\int\limits_{\Gamma}
 \frac{\overline{x-w}}{|x-w|^{n+1}}\,\mathrm{d}\sigma g(x),\,\,\,w\in \Omega,
 \eeq
and
\beq\label{63}
\big(S^{-}[g]\big)(w)=\displaystyle\frac{1}{\bigvee_{n+1}}\int\limits_{\Gamma}
 \frac{\overline{x-w}}{|x-w|^{n+1}}\,\mathrm{d}\sigma g(x)=0,\,\,\,w\in \Omega^{-}=\mathcal{R}^{n+1}\backslash \overline{\Omega}.
 \eeq
Secondly, if (\ref{63}) holds, then by Sochocki-Plemelj formulae (\ref{231}) we know that
(\ref{62}) is just the solution of  Dirichlet  problem (\ref{61}). \vspace{1mm}

\bth
\label{th61}
  If $g\in H(\Gamma)$ then the Dirichlet problem $(\ref{61})$ has the unique solution $(\ref{62})$ if and only uf $(\ref{63})$ is
satisfied.
\eth

By Theorem \ref{T28}, Theorem \ref{th61} can be written in another form as follows.

\bth \label{th62}
 If $g\in H(\Gamma)$ then the Dirichlet  problem $(\ref{61})$ has the unique solution $(\ref{62})$ if and  only if
\beq \label{64}
\big(S[g]\big)(t)=\displaystyle\frac{1}{\bigvee_{n+1}}\int\limits_{\Gamma}
 \frac{\overline{x-t}}{|x-t|^{n+1}}\,\mathrm{d}\sigma g(x)=\frac{1}{2} g(t),\,\,\,\,\,\,t\in\Gamma\vspace{-1mm}
\eeq holds.
\eth

Theorem \ref{th61} and Theorem \ref{th62} can be stated together.

\bth[Conditions for $H$ continuous  boundary value of regular function] \label{th63}
A function $f\!\in\!H(\Gamma)$ is the boundary value of a regular function $\Phi$ in $\Omega$ if and only if $f$ satisfies $(\ref{63})$ or $(\ref{64})$.
\eth

{\bf Dirichlet  problem with  continuous boundary  function.}
Find a left regular function $\Phi$ on $\Omega$,  satisfying the boundary value condition
\beq \label{65}
  \Phi^+(x)=g(x),\,\,\,\, x\in\Gamma,\eeq
where $g\in C(\Gamma)$ is given. \vspace{1mm}

In this case, since  $g \in C(\Gamma),$ the  Sochocki-Plemelj formula is failed in general. To obtain the corresponding result to Theorem \ref{th61} we need to find a substitute of the Sochocki-Plemelj formula. To do so, some tools should be introduced (see \cite{ld}).

\ble[Quasi normal field \mbox{\rm \cite{ld}}] \label{l61}
  Suppose that $\Gamma\subset\mathcal{R}^{n+1}$ is a compact smooth surface and $c\in(0,1)$ is a constant. Then there exists a $C^1$-mapping
  $M\!: \Gamma\rightarrow S^{n}$ such that $\big(M(p),N(p)\big)\geq c$ for each $p\in \Gamma$, where $S^{n}$ is the unit sphere in $\mathcal{R}^{n+1}$, \vspace{1mm}
  $N$ is a continuous unit  on $\Gamma$.
  \ele

  The direction field $M(p)$ $(p\in\Gamma)$ in the above lemma is called a $C^{1}$-unit $c$-quasi normal field on $\Gamma$ associated with the  normal field $N$.

\ble[Tubular neighborhood formed quasi normal lines \mbox{\rm \cite{ld}}] \label{l62}
  Let $\Gamma\subset\mathcal{R}^{n+1}$ be a  compact $C^1$-surface and $M$  a $C^{1}$-unit $c$-quasi normal field on $\Gamma$. Then there exists a number $h>0$
  such that, for any $\epsilon\in(0,h)$,

  $(1)$ $I_{p,\epsilon}\cap I_{q,\epsilon}=\emptyset$, where $p,q\in\Gamma$, $p\not=q$ and
  \beq \label{68}
  I_{p,\epsilon}=\big\{p+\lambda M(p), \,\,-\epsilon<\lambda<\epsilon\big\};
  \eeq

 $(2)$ $\mbox{$\cal V$}_{\epsilon}$ is an open neighborhood  in $\mathcal{R}^{n+1}$ of $\Gamma$, where
 \beq \label{89}
 \mbox{$\cal V$}_{\epsilon}=\big\{p+\lambda M(p),\,\,\, p\in \Gamma,\, -\epsilon<\lambda<\epsilon\big\}.
 \eeq
\ele
\bre\label{R43}{
\rm Clearly, for $\epsilon<h$,
\beq \label{443}
\mbox{$\cal V$}_{\epsilon}={\bigcup}_{p\in \Gamma}I_{p,\epsilon}={\bigcup}_{p\in\Gamma}\big\{p+\lambda M(p), -\epsilon<\lambda<\epsilon\big\}, \,\,\,\,I_{p,\epsilon}\bigcap I_{q,\epsilon}=\emptyset\,\,(p\not=q). \vspace{-2mm}
\eeq
That is to say, through each point of the  neighborhood $\mbox{$\cal V$}_{\epsilon}$ of $\Gamma$
there passes a unique quasi normal line $I_{p,\epsilon}$ of length $2\epsilon$
 with respect to  $C^{1}$-quasi normal $M(p)$. $\mbox{$\cal V$}_{\epsilon}$ is then called a tubular neighborhood of $\Gamma$
  with respect to  $C^{1}$-unit $c$-quasi normal $M(p)$, briefly a tubular neighborhood of
   $\Gamma$. For convenience, we call $h$  an allowable width number of the
 tubular neighborhood $\mbox{$\cal V$}_{\epsilon}$. }
\ere

\ble[Inner and  exterior quasi normal\mbox{\rm\cite{ld}}] \label{l63}
If
 $N(p)$ is the  unit inner normal and $M(p)$ is a $C^{1}$-unit $c$-quasi normal field on $\Gamma$ associated with $N$, then
 \beq
  \big(p,p+\ell M(p)\big)\subseteq\Omega^{+},\,\,\,\,\,
  \big(p,p-\ell M(p)\big)\subseteq\Omega^{-},
  \eeq
 where $\ell$ is some a positive real number, $(a,b)$ stands for
 the line segment from $a$
to $b$.
\ele

Let $\mbox{$\cal V$}_\ell$ be a tubular neighborhood of $\Gamma$. Now we construct a function on
  $\mbox{$\cal V$}^{+}_\ell=\mbox{$\cal V$}_\ell\cap\Omega$ as follows. For $w\in \mbox{$\cal V$}^{+}_\ell$, by the above Lemmas and Remark,
  we know that there exist unique $p\in\Gamma$ and $\lambda\in (0,\ell)$, such that
\beq \label{610}
w=p+\lambda M(p).
\eeq
Then let
 \beq \label{611}
 w_M=p-\lambda M(p),
 \eeq
 which is called the symmetric point of $w$ about $p$ on the line segment $\big(p+h M(p), p-h M(p)\big),$ where $h$ is an allowable width number.
 We have

 \beq\label{612}
 w\in\Omega^{+},\,\,\,\,\,w_M\in \Omega^{-}.
 \eeq
  So, let
  \beq
 \mbox{$\cal V$}_{\ell}^+=\mbox{$\cal V$}_\ell\cap\Omega\,\,\,\, (\ell<h)
  \eeq
that is an open set,  for $f\in C(\Gamma)$ we construct a function $\Delta [f]\!: \mbox{$\cal V$}^{+}_\ell\rightarrow C(V_n)$ by
\beq\label{613}
\Big(\!\Delta\mathcal{S} [f]\!\Big)(w)=\Big(\!S[f]\!\Big)(w)-\Big(\!S[f]\!\Big)\big(w_M\big), \,\,\, \,\,w\in \mbox{$\cal V$}_\ell^+,
\eeq
that is called the symmetric difference of $S$.

In \cite{ld}, we obtain the following result, which is a form of Sochocki-Plemelj formula in case of the continuous density.

\bth \label{th64}
If $f\in C(\Gamma)$, then
\beq \label{864}
\lim_{w\rightarrow p,\, w\in \mbox{\tiny $\cal V$}^{+}_{\ell}}\Big(\!\Delta\mathcal{S}[f]\!\Big)(w)=f(p),\,\,\,\,\,p\in\Gamma,
\eeq
i.e.,
\beq
\lim_{w\rightarrow p,\, w\in \mbox{\tiny $\cal V$}^{+}_{\ell}}\left[\Big(\!S[f]\!\Big)(w)-\Big(\!S[f]\!\Big)\big(w_M\big)\right]=f(p),\,\,\,\,\,p\in\Gamma,
\eeq
and the limit converges uniformly with respect to $p\in \Gamma$.
\eth

By using the above Sochocki-Plemelj formula, now we can give the solution and the conditions of solvability for the Dirichlet  problem  (\ref{65}) with  continuous boundary  function.

\bth
\label{th65} If $f\in C(\Gamma)$, then the
  Dirichlet  problem $(\ref{65})$ has the unique solution $(\ref{62})$ if and only if $(\ref{63})$ is
satisfied.
\eth

 \bth[Conditions for continuous  boundary value of a regular function] \label{th66}
A necessary and sufficient condition for a function $f\!\in\!H(\Gamma)$ to be the boundary value of a function $\Phi$ regular in $\Omega$ is
$(\ref{63})$.
\eth

\vspace{4mm}

\mbox{\bf Acknowledgements}

The first author would like to extend his thankfulness to Professor Qian Tao for his support to the author  during the his stay in University of Macau while the article was prepared.



\vspace{24mm}

\end{document}